%
%
%
%
%
%
%
\magnification=\magstep1
\input amstex
\documentstyle{amsppt}
\voffset-3pc

\define\C{{\Bbb C}}

\define\dee{\partial}
\redefine\O{\Omega}
\redefine\phi{\varphi}
\define\Obar{\overline{\Omega}}
\define\Ot{\widetilde\Omega}
\define\Oh{\widehat\Omega}

\NoRunningHeads
\topmatter
\title
Quadrature domains and kernel function zipping
\endtitle
\author Steven R. Bell${}^*$ \endauthor
\thanks ${}^*$Research supported by NSF grant DMS-0305958 \endthanks
\address
Mathematics Department, Purdue University, West Lafayette, IN  47907 USA
\endaddress
\email bell\@math.purdue.edu \endemail
\abstract
It is proved that quadrature domains are ubiquitous in a very strong
sense in the realm of smoothly bounded multiply connected domains in the
plane.  In fact, they are so dense that one might as well assume that any
given smooth domain one is dealing with is a quadrature domain, and this
allows access to a host of strong conditions on the classical kernel functions
associated to the domain.  Following this string of ideas leads to the
discovery that the Bergman kernel can be ``zipped'' down to a strikingly
small data set.

It is also proved that the kernel functions associated to a quadrature
domain must be algebraic.
\endabstract
\endtopmatter
\document

\hyphenation{bi-hol-o-mor-phic}
\hyphenation{hol-o-mor-phic}

\subhead 1. Introduction \endsubhead
In this paper, we will refine results of B.~Gustafsson in light of
recent results in \cite{8} about the complexity of the classical
kernels functions to show that quadrature domains in the plane are
so dense that one cannot possibly devise a test to determine if a
given smooth domain is a quadrature domain.  The combined methods of
Gustafsson \cite{12} and \cite{8} will also yield a method to
``zip'' the Bergman kernel down to a very small data set consisting
of finitely many complex numbers plus the boundary values of a
{\it single\/} holomorphic function, which I
would venture to christen a {\it Gustafsson function}.  These results
are all a natural outgrowth of the work of Aharonov and Shapiro
\cite{1} and Shapiro \cite{15}, and one consequence of the
Aharonov-Shapiro theorem that Ahlfors maps associated to quadrature
domains are algebraic will be that the Bergman and Szeg\H o kernels
associated to a quadrature domain are algebraic functions.

For the purposes of this paper, we shall call an $n$-connected domain
$\O$ in the plane such that no boundary component is a point a {\it
quadrature domain\/} if there exist finitely many points
$\{w_j\}_{j=1}^N$ in the domain and non-negative integers $n_j$ such
that complex numbers $c_{jk}$ exist satisfying
$$\int_\O f\ dA = \sum_{j=1}^N\sum_{k=0}^{n_j} c_{jk} f^{(k)}(w_j)\tag1.1$$
for every function $f$ in the Bergman space of square integrable
holomorphic functions on $\O$.  Here, $dA$ denotes Lebesgue area
measure.  Many of our results require the function $h(z)\equiv1$ to be
in the Bergman space, and so we shall often also assume that the domain
under study has {\it finite area}.  We remark that there are results
of Sakai that show that, under certain weaker assumptions, a quadrature
domain does have finite area, so some of our results could be stated
with weaker hypotheses.  (See \cite{12} for an explanation of how
Sakai's results relate to the type of quadrature domains we study
here.)

The Ahlfors map associated to a point $a$ in an $n$-connected domain $\O$
such that no boundary component is a point is the holomorphic
function $f_a$  such that $f_a$ maps $\O$ into the unit disc
maximizing the quantity $|f_a'(a)|$ with $f_a'(a)$ real and positive.
This map is an $n$-to-one (counting multiplicities) proper holomorphic
mapping of $\O$ onto the unit disc.

Quadrature domains have particularly simple kernel functions, as our
first theorem shows.

\proclaim{Theorem 1.1}
Suppose that $\O$ is an $n$-connected quadrature domain of finite area
in the plane such that no boundary component is a point.  Then the
Bergman kernel function $K(z,w)$ associated to $\O$ is a rational
combination of two Ahlfors maps $f_a$ and $f_b$ in the sense that
$K(z,w)$ is a rational combination of $f_a(z)$, $f_b(z)$,
$\overline{f_a(w)}$, and $\overline{f_b(w)}$.
The same can be said of the square $S(z,w)^2$ of the Szeg\H o kernel.
Furthermore, the classical functions $F_j'$ are rational functions of
$f_a$ and $f_b$.
\endproclaim

The functions $F_j'$ are defined precisely in \S2.

Aharonov and Shapiro \cite{1} proved that Ahlfors maps associated to
quadrature domains are algebraic.  Hence, Theorem~1.1 yields that
quadrature domains also have algebraic kernel functions.

\proclaim{Theorem 1.2}
Suppose that $\O$ is an $n$-connected quadrature domain in the plane
of finite area such that no boundary component is a point. The Bergman
and Szeg\H o kernel functions associated to $\O$ are algebraic functions.
The classical functions $F_j'$ are also algebraic.
\endproclaim

Similar statements to Theorems~1.1 and~1.2 hold for the Poisson kernel
and first derivative of the Green's function.  These results follow
from formulas appearing in \cite{6} and we do not spell them out
here.

Ahlfors maps extend to the double (as described in \S2 of this paper),
and it follows that, under the hypotheses of Theorem~1.1, the Bergman
kernel function extends meromorphically to the double $\Oh$ of $\O$,
and is therefore a rational combination of any two functions that
generate the meromorphic functions on the double, i.e., any two
functions that form a primitive pair for the double.  The Bergman
kernel always extends to the double as a meromorphic differential,
but extending as a meromorphic functions is rather unusual behavior
for the kernel.  This condition leads to a number of other strong
conclusions that we now begin to enumerate.

\proclaim{Theorem 1.3}
Suppose that $\O$ is an $n$-connected quadrature domain of finite area
in the plane such that no boundary component is a point.  If $f$ is any
proper holomorphic mapping of $\O$ onto the unit disc, then $f'$ extends to
the double of $\O$ as a meromorphic function.
\endproclaim

Under the assumptions of Theorem~1.3, since both $f$ and $f'$ extend
to the double, they are algebraically dependent, i.e., there exists
an irreducible polynomial $P(z,w)$ on $\C^2$ such that $P(f',f)\equiv0$
on $\O$.  This was proved by other means by Aharonov and Shapiro in
\cite{1}.  It is proved in \cite{8} that the condition $P(f',f)=0$
has a number of implications, two of which are that the kernel functions
are generated by only two functions and that the kernel functions
extend to a compact Riemann surface.  In the setting of Theorem~1.3,
however, we have the stronger conclusion that the Bergman kernel extends to
the compact Riemann surface which is the double of $\O$, and that it
is generated by any two functions that form a primitive pair for the
double.

When combined with the main theorem of \cite{9}, Theorem~1.3 yields
that the infinitesimal Carath\'eodory metric associated to an
$n$-connected quadrature domain of finite area such that no boundary
component is a point is a real algebraic function which is a rational
combination of two Ahlfors maps and their conjugates.

It is shown in \cite{7} that, under the assumptions of Theorem~1.3,
if $f$ is any proper holomorphic map of $\O$ onto the unit disc, then
it is possible to find an Ahlfors map $f_b$ such that $f$ and $f_b$
extend to the double and generate the meromorphic functions on the
double.  Hence, it follows that $f'=R(f,f_b)$ for some rational
function.  Also, we may conclude that, given a proper map $f$, the
Bergman kernel is a rational combination of $f$ and some other Ahlfors
map.  Furthermore, since both $f'$ and $f_b'$ extend meromorphically
to the double, we deduce the next rather odd sounding theorem.

\proclaim{Theorem 1.4}
Suppose that $\O$ is an $n$-connected quadrature domain of finite area
in the plane such that no boundary component is a point.  If $H$ is any
meromorphic function on $\O$ that extends meromorphically to the
double of $\O$, then $H'$ also extends meromorphically to the double
of $\O$.  Furthermore, $H$ is algebraic.
\endproclaim

Of course if $H$ is meromorphic on $\O$, then $H'$ is meromorphic
on $\O$.  The content of the theorem is that $H'$ {\it extends\/}
meromorphically to the double if $H$ does.

The property in Theorem~1.4 turns out to characterize a class of
generalized quadrature domains, and we explore this line of
reasoning in \S5.

When we combine the ideas used in the proofs of the results
above with those of Aharonov and Shapiro \cite{1}
and Gustafsson \cite{12}, we can show that the kernel functions
associated to quadrature domains are particularly simple when
restricted to the boundary.

\proclaim{Theorem 1.5}
Suppose that $\O$ is an $n$-connected quadrature domain in the plane
of finite area such that no boundary component is a point. The Bergman
kernel $K(z,w)$ and the square $S(z,w)^2$ of the Szeg\H o kernel are
rational functions of $z$, $\bar z$, $w$, and $\bar w$ on
$b\O\times b\O$ minus the boundary diagonal.  The functions $F_j(z)$
are rational functions of $z$ and $\bar z$ when restricted to the
boundary.  Furthermore, the unit tangent vector function $T(z)$ is
such that $T(z)^2$ is a rational function of $z$ and $\bar z$ for
$z\in b\O$.
\endproclaim

I had conjectured in \cite{5} that every $n$-connected domain in
the plane such that no boundary component is a point is conformally
equivalent to a domain with algebraic kernel functions.  Jeong and
Taniguchi \cite{14} recently verified this conjecture.  Since
Gustafsson proved in \cite{12} that every such domain is conformally
equivalent to a quadrature domain of finite area, Theorem~1.2
gives an alternate way of seeing that every $n$-connected domain in
the plane such that no boundary component is a point is conformally
equivalent to a domain with algebraic kernel functions.

We remark here that, although the last part of Theorem~1.5 might seem
to suggest that the Bergman kernel associated to a quadrature domain
could be a simple rational function of some kind, it can never happen
that the Bergman kernel is a rational function in the setting of
multiply connected domains (see \cite{4}).

The main results of this paper together with Gustafsson's theorem
that any finitely connected domain in the plane such that no boundary
component is a point is conformally equivalent to a smoothly bounded
quadrature domain suggests that quadrature domains might serve to
play a role in the multiply connected setting similar to that played
by the unit disc for simply connected regions.

Quadrature domains with smooth boundaries are particularly
appealing and we can refine arguments of Gustafsson to prove
the next theorem, which shows that very fine modifications
can be made to any smoothly bounded domain to make it a quadrature
domain.

\proclaim{Theorem 1.6}
Suppose that $\O$ is a bounded $n$-connected domain whose
boundary consists of $n$ non-intersecting $C^\infty$ smooth simple
closed curves.  There is a meromorphic function $g$ on the double
of $\O$ which has no poles on $\Obar$ such that $g$ is as close to
the identity map in $C^\infty(\Obar)$ as desired.  The domain
given by $g(\O)$ is a quadrature domain which is $C^\infty$ close
to $\O$ and conformally equivalent to $\O$.
\endproclaim

The analytic objects attached to the quadrature domain $g(\O)$ have
the strong extension properties given in the preceding theorems
and they are $C^\infty$ close to the analytic objects attached to
$\O$.  In particular, Aharonov and Shapiro \cite{1} (with some
refinements by Gustafsson \cite{12}) showed that the boundary of
$g(\O)$ is an algebraic curve minus perhaps finitely many points.
Thus, the proof of Theorem~1.6 will yield a concrete method to
approximate in $C^\infty$ a non-intersecting group of $n$ simple
closed $C^\infty$ curves by an algebraic curve.  In fact, the algebraic
curve can be described by $|f(z)|^2=1$ where $f$ is any Ahlfors map.

We describe in \S4 how the Bergman kernel can be recovered from
the boundary values of $g$ in a very simple and efficient manner.

Gustafsson proved that the function $g$ in Theorem~1.6 maps $\O$ to a
quadrature domain.  We shall show that $g(z)$ can be taken to be a
linear combination of functions of the form $S(z,b)/L(z,a)$ where
$S(z,b)$ is the Szeg\H o kernel and $L(z,a)$ is the Garabedian
kernel, and $b$ ranges over a small open subset of $\O$ while $a$
is fixed.  Consequently, we shall be able to restrict the points
$w_j$ in the defining property (1.1) of the quadrature domain
$g(\O)$ to a
small set.   We shall also be able to specify the numbers $n_j$ in
rather surprising ways.  In particular, we shall be able to stipulate
that $n_j=1$ for each $j$.  Thus, any smooth domain is conformally
equivalent to a nearby quadrature domain where the simple point
masses are contained in an arbitrarily small arbitrary disc that
is compactly contained in $\O$.  Another way to state this is to say
that it is possible to strongly approximate the two dimensional field
generated by a uniform charge density on a smoothly bounded plate
with holes by point charges at finitely many points in an arbitrarily
small open subset of the plate.  This result is stated precisely in
the following theorem.

\proclaim{Theorem~1.7}
Suppose that $\O$ is a bounded $n$-connected domain whose
boundary consists of $n$ non-intersecting $C^\infty$ smooth simple
closed curves.  Let $D_\epsilon(w_0)$ be any disc which is compactly
contained in $\O$.  There is a quadrature domain which is $C^\infty$ close
to $\O$ and conformally equivalent to $\O$ such that the point masses
appearing in (1.1) all fall in $D_\epsilon(w_0)$ and have weight $n_j=1$.
Furthermore, given $w_0$ in $\O$, there is a quadrature domain which
is $C^\infty$ close to $\O$ and conformally equivalent to $\O$ such
that $w_0$ is the only point mass appearing in (1.1), i.e., $N=1$ in
(1.1).
\endproclaim

In \S6 of this paper, we show how many of the same ideas can
be extended to quadrature domains with respect to boundary arc
length measure.

\subhead 2. Preliminaries\endsubhead
It is a standard construction in the theory of conformal mapping to
show that an $n$-connected domain $\O$ in the plane such that no boundary
component is a point is conformally equivalent via a map $\Phi$ to a
bounded domain $\Ot$ whose boundary consists of $n$ simple closed
$C^\infty$ smooth real analytic curves.  Since such
a domain $\Ot$ is a bordered Riemann surface, the double of $\Ot$
is an easily realized compact Riemann surface.  We shall say that
an analytic or meromorphic function $h$ on $\O$ {\it extends meromorphically
to the double of\/} $\O$ if $h\circ\Phi^{-1}$ extends meromorphically
to the double of $\Ot$.  Notice that whenever $\O$ is itself a bordered
Riemann surface, this notion is the same as the notion that $h$
extends meromorphically to the double of $\O$.  We shall say that
two functions $G_1$ and $G_2$ extend to the double and generate the
meromorphic functions on the double of $\O$, and that they therefore
form a primitive pair for the double of $\O$, if $G_1\circ\Phi^{-1}$ and
$G_2\circ\Phi^{-1}$ extend to the double of $\Ot$ and form a primitive
pair for the double of $\Ot$ (see Farkas and Kra \cite{11} for the
definition and basic facts about primitive pairs).

It is proved in \cite{7} that if $\O$ is an $n$-connected domain
in the plane such that no boundary component is a point, then almost
any two distinct Ahlfors maps $f_a$ and $f_b$ generate the meromorphic
functions on the double of $\O$.  It is also proved that any proper
holomorphic mapping from $\O$ to the unit disc extends to the double
of $\O$.

Suppose that $\O$ is a bounded $n$-connected domain whose
boundary consists of $n$ non-intersecting $C^\infty$ smooth simple
closed curves.  The Bergman kernel $K(z,w)$ associated to $\O$ is
related to the Szeg\H o kernel via the identity
$$K(z,w)=4\pi S(z,w)^2+\sum_{i,j=1}^{n-1}
A_{ij}F_i'(z)\overline{F_j'(w)},\tag 2.1$$
where the functions $F_i'(z)$ are well known classical functions of
potential theory described as follows.  The harmonic function $\omega_j$
which solves the Dirichlet problem on $\O$ with boundary data equal to
one on the boundary curve $\gamma_j$ and zero on $\gamma_k$ if $k\ne j$
has a multivalued harmonic conjugate.  Let $\gamma_n$ denote the outer
boundary curve.  The function $F_j'(z)$ is a single valued holomorphic
function on $\O$ which is locally defined as the derivative of
$\omega_j+iv$ where $v$ is a local harmonic conjugate for $\omega_j$.
The Cauchy-Riemann equations reveal that $F_j'(z)=2(\dee\omega_j/\dee z)$.

The Bergman and Szeg\H o kernels are holomorphic in the first variable
and antiholomorphic in the second on $\O\times\O$ and they are hermitian,
i.e.,  $K(w,z)=\overline{K(z,w)}$.  Furthermore, the Bergman and
Szeg\H o kernels are in
$C^\infty((\Obar\times\Obar)-\{(z,z):z\in b\O\})$ as functions of $(z,w)$
(see \cite{2, page~100}).

We shall also need to use the
Garabedian kernel $L(z,w)$, which is related to the Szeg\H o
kernel via the identity
$$\frac{1}{i} L(z,a)T(z)=S(a,z)\qquad\text{for $z\in b\O$ and $a\in\O$}
\tag2.2$$
where $T(z)$ represents the complex unit tangent vector at $z$ pointing
in the direction of the standard orientation of $b\O$.
For fixed $a\in\O$, the kernel $L(z,a)$ is a holomorphic function of $z$
on $\O-\{a\}$ with a simple pole at $a$ with residue $1/(2\pi)$.
Furthermore, as a function of $z$, $L(z,a)$ extends to the boundary
and is in the space $C^\infty(\Obar-\{a\})$.  In fact, $L(z,w)$
is in $C^\infty((\Obar\times\Obar)-\{(z,z):z\in\Obar\})$ as a function
of $(z,w)$ (see \cite{2, page~102}).  Also, $L(z,a)$ is non-zero for all
$(z,a)$ in $\Obar\times\O$ with $z\ne a$ and $L(a,z)=-L(z,a)$ (see
\cite{2, page~49}).

For each point $a\in\O$, the function of $z$ given by
$S(z,a)$ has exactly $(n-1)$ zeroes in $\O$ (counting multiplicities) and
does not vanish at any points $z$ in the boundary of $\O$ (see
\cite{2, page~49}).

Given a point $a\in\O$, the Ahlfors map $f_a$ associated to the pair $(\O,a)$
is a proper holomorphic mapping of $\O$ onto the unit disc.  It is an
$n$-to-one mapping (counting multiplicities), it extends to be in
$C^\infty(\Obar)$, and it maps each boundary curve $\gamma_j$ one-to-one
onto the unit circle.  Furthermore, $f_a(a)=0$, and $f_a$ is the unique
function mapping $\O$ into the unit disc maximizing the quantity $|f_a'(a)|$
with $f_a'(a)>0$.  The Ahlfors map is related to the Szeg\H o kernel
and Garabedian kernel via (see \cite{2, page~49})
$$f_a(z)=\frac{S(z,a)}{L(z,a)}.\tag2.3$$
Note that $f_a'(a)=2\pi S(a,a)\ne 0$.  Because $f_a$ is $n$-to-one, $f_a$
has $n$ zeroes.  The simple pole of $L(z,a)$ at $a$ accounts for the simple
zero of $f_a$ at $a$.   The other $n-1$ zeroes of $f_a$ are given by the
$(n-1)$ zeroes of $S(z,a)$ in $\O-\{a\}$.

When $\O$ does not have smooth boundary, we define the kernels and domain
functions above as in \cite{6} via a conformal mapping to a domain with
real analytic boundary curves.

\subhead 3. Proofs of the theorems\endsubhead
If $\O$ is an $n$-connected quadrature domain of finite area in the plane
such that no boundary component is a point, then the Bergman kernel
function associated to $\O$ satisfies an identity of the form
$$1\equiv\sum_{j=1}^N\sum_{m=0}^{n_j} c_{jm} K^{(m)}(z,w_j)
\tag3.1$$
where $K^{(m)}(z,w)$ denotes $(\dee^m/\dee\bar w^m)K(z,w)$ and
the points $w_j$ are the points that appear in the characterizing
formula (1.1) of quadrature domains.  This observation is
usually attributed to Avci in his unpublished Stanford PhD thesis.
It can be seen by noting that the inner product of an analytic
function against the function $h(z)\equiv1$ and against the sum on
the right hand side of (3.1) agree for all functions in the Bergman
space.  Hence the two functions must be equal.  Note that we must
assume that $\O$ has finite area here just so that $h(z)\equiv1$ is
in the Bergman space.

\demo{Proof of Theorem~1.1}
Since the Bergman kernel is equal to $(\dee^2/\dee z\dee\bar w)$
of the Green's function, functions that are of the form of the right
hand side of (3.1) belong to the the class $\Cal A$ of
\cite{8, p.~20}.  Hence, the function $A(z)\equiv1$ belongs to
$\Cal A$.  Theorem~2.3 of \cite{8} states that if
$G_1$ and $G_2$ are any two meromorphic functions on $\O$ that extend
to the double of $\O$ to form a primitive pair and if $A(z)$ is any
function from the class~$\Cal A$ other than the zero function, then
the Bergman kernel associated to $\O$ can be expressed as
$$K(z,w)=A(z)\overline{A(w)}
R_1(G_1(z),G_2(z),\overline{G_1(w)},\overline{G_2(w)})$$
where $R_1$ is a complex rational function of four complex variables.
Similarly, the Szeg\H o kernel can be expressed as
$$S(z,w)^2= A(z)\overline{A(w)}
R_2(G_1(z),G_2(z),\overline{G_1(w)},\overline{G_2(w)})$$
where $R_2$ is rational, and
the functions $F_j'$ can be expressed
$$F_j'(z)=A(z)R_3(G_1(z),G_2(z))$$
where $R_3$ is rational.  Furthermore, every proper holomorphic
mapping of $\O$ onto the unit disc is a rational combination of
$G_1$ and $G_2$.  It therefore follows now that the Bergman kernel
is a rational combination of any two meromorphic functions on
$\O$ that extend to the double to form a primitive pair.  Since
almost any two distinct Ahlfors maps form a primitive pair (see
\cite{7}), the proof of Theorem~1.1 is complete.
\enddemo

\demo{Proof of Theorem~1.3}
Suppose that $\O$ is an $n$-connected quadrature domain in the plane
such that no boundary component is a point and suppose that $f$ is a
proper holomorphic mapping of $\O$ onto the unit disc.  We may
compose $f$ with a M\"obius transformation $\phi$ so that $F=\phi\circ f$
has only simple zeroes at, say $a_1,a_2,\dots,a_N$, where $N$ is the
order of the proper map~$f$.  It is proved in \cite{2, p.~65} that
the Bergman kernel transforms under this proper map according to
$$F'(z)K_D(F(z),0)=\sum_{n=1}^N K(z,a_n)/\/\overline{F'(a_n)}$$
where $K_D(z,w)=\pi^{-1}(1-z\bar w)^{-2}$ is the Bergman kernel
for the unit disc.  Notice that $K_D(z,0)\equiv \pi^{-1}$, and
so it follows that $F'(z)$ is given by a linear combination of
functions of the form $K(z,a_n)$, and thus $F'$ extends to the
double of $\O$ by Theorem~1.1.  But $F'(z)=f'(z)\phi'(f(z))$,
and since $\phi$ is rational and $f$ extends to the double of
$\O$, it now follows that $f'(z)$ extends to the double of $\O$.
The proof is complete.
\enddemo

\demo{Proof of Theorem~1.5}
In the setting of Theorem~1.5, Gustafsson \cite{12} generalized
a result of Aharonov and Shapiro \cite{1} to prove that the boundary of
$\O$ is given by an algebraic curve and that there exists a function
$H(z)$ which is meromorphic on $\O$ with continuous boundary values
such that $H(z)=\bar z$ on $b\O$.  Let $G(z)=z$.  Gustafsson proved
that $H(z)$ and $G(z)$ extend to the double of $\O$ to form a
primitive pair.  Hence, there exists an irreducible polynomial
$P(z,w)$ on $\C^2$ such that $P(H(z),G(z))\equiv0$ on $\O$.  This
shows that $H(z)$ is an algebraic function of $z$.  We know
that the Bergman kernel is generated by $z$ and $H(z)$.  Hence,
this gives another way to see that the Bergman kernel is algebraic.
It is proved in \cite{5} that if the Bergman kernel is algebraic,
then so is the Szeg\H o kernel, all proper holomorphic maps onto
the unit disc, and the classical functions $F_j'$.

Now since the kernels $K(z,w)$ and $S(z,w)^2$ and the proper
holomorphic maps to the unit disc and the functions $F_j$
are all generated by $G(z)$ and $H(z)$, and since these functions
are equal to $z$ and $\bar z$, respectively on the boundary, we may
deduce most of the rest of the claims made in Theorem~1.5.  To
finish the proof, note that identity (2.2) yields that
$$T(z)^2=-\frac{S(a,z)^2}{L(z,a)^2}$$
where $a$ is an arbitrary point chosen and fixed in $\O$.
The function $S(z,a)^2$ is a rational function of $z$ and
$\bar z$ on the boundary.  Identity (2.3) yields that
$L(z,a)^2=S(z,a)^2/f_a(z)^2$, and so $L(z,a)^2$ is also a rational
function of $z$ and $\bar z$ on the boundary.  Finally, it follows
that $T(z)^2$ is a rational function of $z$ and $\bar z$.
\enddemo

We remark that, since the antiholomorphic Schwarz reflection function
$S(z)$ across a real analytic boundary curve of $\O$ satisfies
$\overline{f(S(z))}=1/f(z)$ when $f$ is a proper holomorphic mapping
onto the unit disc, it follows that $S(z)$ is algebraic whenever proper
holomorphic maps to the disc are.

\demo{Proof of Theorem~1.6}
Suppose that $\O$ is a bounded $n$-connected domain whose boundary
consists of $n$ non-intersecting $C^\infty$ smooth simple closed
curves.  I proved in \cite{10} (see also \cite{2, p.~29}) that the
complex linear span
of the set of functions of $z$ given by $\{S(z,b):b\in\O\}$ is dense in
$A^\infty(\O)$, the subset of $C^\infty(\Obar)$ consisting of holomorphic
functions on $\O$.  The proof given there is constructive.
It is proved in \cite{3} that there is a dense open set of points
$a$ in $\O$ such that $S(z,a)$ has $n-1$ simple zeroes as a function of $z$.
Fix such a point $a$ and let $a_1,a_2,\dots,a_{n-1}$ denote the zeroes
of $S(z,a)$.  The functions of $z$ given by $S(z,b)/L(z,a)$ extend
meromorphically to the double $\Oh$ of $\O$ because identity (2.2)
shows that $S(z,b)/L(z,a)$ agrees with the conjugate of $L(z,b)/S(z,a)$
on the boundary of $\O$.  Let $R(z)$ denote the antiholomorphic
reflection function which maps $\O$ to its reflected copy in the double.
Notice that the extended function has no poles in $\Obar$ and, if $b$ is
not equal to any of the zeroes $a_j$, then it has
simple poles at $R(b)$ and the points $\{R(a_j)\}_{j=1}^{n-1}$.

The function  $H(z)$ which is equal to $(z-a)L(z,a)$ 
for $z\in\O$, $z\ne a$, and equal to $1/2\pi$ at $z=a$ is in
$A^\infty(\O)$.  Hence, we may find finitely many points $b_j$ in $\O$ such
that a linear combination  $\sum_{j=1}^N c_{j}S(z,b_j)$ is as close to
$H(z)$ in $A^\infty(\O)$ as desired.  Now the function $g(z)$ given by
$$a+\sum_{j=1}^N c_{j} S(z,b_j)/L(z,a)$$
extends to be a meromorphic function on the double of $\O$ which is
close in $C^\infty(\Obar)$ to the identity function.  It is this
function $g$ that we wish to call a {\it Gustafsson function}.
We shall use it in the next section to zip the Bergman kernel.

Gustafsson proved in \cite{12} that the poles of the function
$g$ on the reflected copy of $\O$ in the double of $\O$ reflect
back to the points in $\O$ that map under $g$ to points that appear
in the quadrature identity for $g(\O)$.  Hence, the points in $g(\O)$
that would appear in the quadrature identity (1.1) for $g(\O)$ are
among the images under $g$ of the points $a_1,\dots,a_{n-1}$
and $b_1,\dots,b_N$.  We shall refine the proof above to get more
control over these points momentarily.
\enddemo

\demo{Proof of Theorem~1.7}
The proof just given of Theorem~1.6 can be altered so that the points
$b_j$ fall in a very small set and so that all the integers $n_j$ in the
quadrature identity for $g(\O)$ are equal to one.  Indeed, let
$D_\epsilon(w_0)$ be any disc which is compactly contained in $\O$.
Choose $a$ in $\O$ such that $S(z,a)$ has $n-1$ simple zeroes as a
function of $z$.  Let $a_1,a_2,\dots,a_{n-1}$ denote the zeroes of
$S(z,a)$ and let $a_0=a$.
We shall now repeat the argument above, but we shall restrict the
points $b$ to be in $D_\epsilon(w_0)$.  Indeed, we now claim that
the complex linear
span $\Cal L$ of $\{S(z,b)\,:\, b\in D_\epsilon(w_0)\}$ is dense in
$A^\infty(\O)$.  The dual space $A^{-\infty}(\O)$ of $A^\infty(\O)$
is described in \cite{2, p.~117} (see also \cite{10}).  If $\Cal L$
were not dense in $A^\infty(\O)$, then there would be a function $h\in
A^{-\infty}(\O)$ which is not the zero function, but which is orthogonal
to $\Cal L$ with respect to the non-degenerate pairing which extends
the usual $L^2$ inner product on $\O$.  But the function $H(b)$ given
as $\langle h,S(\cdot,b)\rangle$ is a holomorphic function of $b$ on
$\O$.  Thus, if $h$ is orthogonal to $\Cal L$, then $H(b)$ vanishes
on $D_\epsilon(w_0)$, and is therefore zero on all of $\O$.  This
shows that $h$ is orthogonal to $S(z,b)$ for all $b$ in $\O$, and
we know that these functions span a dense subset of $A^\infty(\O)$.
Hence, $h\equiv0$, and this contradiction yields that $\Cal L$ must
be dense.

The function $H(z)$ given by
$$H(z)=zf_a(z)L(z,a)$$
for $z\in\O$, $z\ne a$, and equal to $af_a'(a)/2\pi$ at $z=a$ is in
$A^\infty(\O)$.  Hence, we may find finitely many points $b_j$ in
$D_\epsilon(w_0)$ such that a linear
combination  $L(z)=\sum_{j=1}^N c_{j}S(z,b_j)$ is as close to
$H(z)$ in $A^\infty(\O)$ as desired.  Now the function $g(z)$ given by
$$f_a(z)^{-1}\sum_{j=1}^N c_{j} S(z,b_j)/L(z,a)$$
extends to be a meromorphic function on the double of $\O$ which is
$C^\infty$ close to the identity function near and up to the boundary of
$\O$.  We shall now make some adjustments to this function to eliminate
any poles that might occur at the zeros of $f_a$.  Since the complex
span $\{S(z,b)\,:\, b\in D_\epsilon(w_0)\}$ is dense in
$A^\infty(\O)$, there exist points $B_k$ in $D_\epsilon(w_0)$ such
that $\det [M_{jk}]\ne 0$ where $[M_{jk}]$ is the $n\times n$ matrix
given by  $M_{jk}= S(a_j,B_k)$  in which the indices range over
$j=0,1,\dots n-1$ and $k=0,1,\dots,n-1$.  Since $H(z)$ vanishes
at the zeroes $a_j$, $j=0,1,\dots n-1$, of $f_a$, the
complex numbers $L(a_j)$ are small, and the closer $L(z)$ is to $H(z)$
in $A^\infty(\O)$, the smaller they are.  Let $\mu_{jk}$ solve the
system
$$L(a_j)=\sum_{k=0}^{n-1} \mu_{jk}S(a_j,B_k),$$
for $j=0,1,\dots n-1$.  Note that the complex numbers $\mu_{jk}$ are
small and that they go to zero as $L$ tends to $H$ in $A^\infty(\O)$.
We now revise the definition of the function $g(z)$ to be
$$f_a(z)^{-1}\left(\sum_{j=1}^N c_{j} S(z,b_j)/L(z,a)-
\sum_{k=0}^{n-1} \mu_{jk}S(z,B_k)/L(z,a)\right).$$
This function has the virtue that it has no poles at the zeroes of
$f_a$, and because it is $C^\infty$  close to the identity near the
boundary of $\O$, it is close to the identity in $C^\infty(\Obar)$.
Furthermore, the extension of this function to the reflected side in
the double is given by the conjugate of
$$f_a(z)\left(\sum_{j=1}^N c_{j} L(z,b_j)/S(z,a)-
\sum_{k=0}^{n-1} \mu_{jk}L(z,B_k)/S(z,a)\right),$$
(where we are thinking $z=R(\zeta)$ where $R$ is the reflection
function on the double).
This function has only simple poles at the points $b_j$ and $B_k$
in $D_\epsilon(w_0)$.  This completes the first part of the proof
of Theorem~1.7.  To prove the last assertion in the statement of
Theorem~1.7, repeat the argument above, noting that the same
reasoning shows that the complex linear span of
$\{S^{(m)}(z,w_0), m=0,2,\dots\}$ where
$S^{(m)}(z,w)=(\dee^m/\dee\bar w^m)S(z,w)$ is also dense in $A^\infty(\O)$,
and also observing that identity (2.2) can be used in the same way to
show that $S^{(m)}(z,w_0)/L(z,a)$ extends meromorphically to the double.
\enddemo

\subhead 4. How to zip the Bergman kernel\endsubhead
Suppose that $\O$ is a bounded $n$-connected domain whose
boundary consists of $n$ non-intersecting $C^\infty$ smooth simple
closed curves.  Let $g(z)$ denote a Gustafsson function as
constructed in the proofs of Theorems~1.6 or~1.7.  Let $\Oh$ denote
the double of $\O$ and let $R(z)$ denote the antiholomorphic
reflection function which maps $\O$ to its reflected copy $\Ot$.
Let $G(z)$ denote the meromorphic extension of $g(z)$ to the double.
Gustafsson \cite{12} proved that $G(z)$ and $\overline{G(R(z))}$ form
a primitive pair for the field of meromorphic functions on $\Oh$.
Now $g(\O)$ is a quadrature domain and the function $g(z)$ transforms
to be the function $z$ on $g(\O)$.  Hence, $z$ and
$\overline{G(R(g^{-1}(z)))}$ extend meromorphically to the double of
$g(\O)$ and form a primitive pair.  Let $h(z)$ denote the meromorphic
function $\overline{G(R(g^{-1}(z)))}$.  Notice that $h(z)$ is equal
to $\bar z$ on $b\O$ and that $h$ extends $C^\infty$ smoothly up to
the boundary.

Let $\{w_j\}_{j=1}^N$ denote the finitely many poles of $h(z)$ in $\O$
and let $n_j$ be equal to the order of the pole at $w_j$.  The
numbers $N$, $w_j$, and $n_j$ are exactly the numbers that appear
in (1.1) in the quadrature identity for $g(\O)$.  Let $P_j(z)$ denote
the principal part of $h(z)$ at $w_j$.  Theorem~1.1 yields that the
Bergman kernel associated to $g(\O)$ is a rational combination of $z$
and $h(z)$.  We don't need to zip the function $z$.
Recall that $h(z)=\bar z$ on the boundary of $g(\O)$, and so the
function $h(z)$ can be zipped via the formula
$$h(z)-\sum_{j=1}^N P_j(z)=\frac{1}{2\pi i}\int_{b\O}\frac{\bar
\zeta-\sum_{j=1}^N P_j(\zeta)}{\zeta-z}\ d\zeta.$$
But $\int_{b\O}\frac{1}{(\zeta-w_j)^k(\zeta-z)}\ d\zeta$
is zero for positive integers $k$.  Hence
$$h(z)=\sum_{j=1}^N P_j(z)+\frac{1}{2\pi i}\int_{b\O}\frac{\bar
\zeta}{\zeta-z}\ d\zeta.$$
Define $\Cal Q$ via
$$\Cal Q(z)=\frac{1}{2\pi i}\int_{b\O}\frac{\bar\zeta}{\zeta-z}\ d\zeta.
\tag4.1$$
We may state that $\Cal Q$ is a holomorphic function on $\O$ which extends
meromorphically to the double of $\O$ without poles in $\Obar$ and
that $\Cal Q$ is an algebraic function.  We have just proved that the
Bergman kernel associated to $g(\O)$ is a rational combination of
$z$, $\Cal Q(z)$, $\bar w$, and $\overline{\Cal Q(w)}$.
A rational function is encoded by finitely many complex
coefficients and a few positive integers.  These numbers carry all
the information that is needed to unzip the Bergman kernel for
the quadrature domain $g(\O)$ via formula (4.1).  (Gustafsson \cite{12}
proved that any quadrature domain of finite area can be expressed as
$g(\O)$ for some such $g$ and smooth $\O$, and so this result can be
easily generalized.)

We now turn to zipping the Bergman kernel $K(z,w)$ for $\O$.  Let
$H(z)$ denote the function $\overline{G(R(z))}$, which is meromorphic
on $\O$ and extends $C^\infty$ smoothly up to $b\O$ and has boundary
values equal to $\overline{g(z)}$.  The
transformation formula for the Bergman kernel under biholomorphic
maps together with the form of the Bergman kernel for $g(\O)$ reveals
that $K(z,w)$ is equal to $g'(z)\overline{g'(w)}$ times a rational
function of $g(z)$, $H(z)$,
$\overline{g(w)}$, and $\overline{H(w)}$.  The Cauchy integral formula
$$g'(w)=\frac{1}{2\pi i}\int_{b\O}\frac{g(z)}{(z-w)^2}\ dz$$
allows us to obtain $g'$ inside $\O$ from the boundary values of $g$,
and of course $g$ can be unzipped in the same manner.  The function
$H(z)=\overline{G(R(z))}$ can be recovered in a way to similar to how
we handled $h$ above.  Let $P(z)$ denote the sum of the principal
parts of $H(z)$.  Then, as above,
$$H(z)=P(z)+\frac{1}{2\pi i}\int_{b\O}\frac{\overline{g(\zeta)}}
{\zeta-z}\ d\zeta.$$
Hence we see that the Bergman kernel can be recovered from the
boundary values of the single function $g$, assuming that finitely many
coefficients from two rational functions are known.

\subhead 5. Generalized quadrature domains\endsubhead
We shall call an $n$-connected domain
$\O$ in the plane such that no boundary component is a point a {\it
generalized quadrature domain\/} if there exist finitely many points
$\{w_j\}_{j=1}^N$ in the domain, non-negative integers $n_j$, and
finitely many continuous closed curves or curve segments $\sigma_m$
in $\O$ such that complex numbers $c_{jk}$ and $b_m$ exist satisfying
$$\int_\O f\ dA = \sum_{j=1}^N\sum_{k=0}^{n_j} c_{jk} f^{(k)}(w_j)
+\sum_{m=1}^M b_m \int_{\sigma_m} f(z)\ dz
\tag5.1$$
for every function $f$ in the Bergman space of square integrable
holomorphic functions on $\O$.  Here, $dA$ denotes Lebesgue area
measure.  As before, we shall also need to assume that the domain
under study has {\it finite area}.  The property of being a
generalized quadrature domain and the conditions mentioned in
Theorems~1.1-1.4 are tied together nicely in the following theorem.

\proclaim{Theorem 5.1}
Suppose that $\O$ is an $n$-connected domain in the plane
of finite area such that no boundary component is a point.
The following conditions are equivalent.
{\roster
\item
$\O$ is a generalized quadrature domain.
\item
The Bergman kernel extends to the double of $\O$ as a meromorphic
function, i.e., the Bergman kernel is generated by the restriction
of two functions of one variable that form a primitive pair for the
field of meromorphic functions on the double of $\O$.
\item
There exists a proper holomorphic mapping $f$ of $\O$ onto the unit
disc such that $f'$ extends to the double of $\O$ as a meromorphic
function.
\item
The derivative of every proper holomorphic mapping of $\O$ onto
the unit disc extends to the double of $\O$ as a meromorphic
function.
\item
Every function $H$ on $\O$ that extends meromorphically to the double
of $\O$ is such that $H'$ also extends to the double of $\O$.
\endroster}
\endproclaim

We have proved most of the equivalences in Theorem~5.1 in the proofs
of Theorems~1.1-1.4.  To finish the proof, we need only show that if
$f$ is a proper holomorphic mapping
of $\O$ onto the unit disc such that $f'$ extends meromorphically to
the double of $\O$, then $\O$ is a generalized quadrature domain.  The
condition that $f'$ extends to the double means that there is a
conformal map $\Phi$ from $\O$ to a bounded domain $\Ot$ whose boundary
consists of $n$ simple closed real analytic curves, and
$f'\circ\Phi^{-1}$ extends to the double of $\Ot$.  Let $\phi$
denote the inverse of $\Phi$.  Since $f\circ\phi$ is a proper
holomorphic mapping of $\Ot$ onto the unit disc, and since these two
domains have real analytic boundary, the mapping $f\circ\phi$ extends
holomorphically past the boundary of $\Ot$.  Hence, the derivative
$\phi'\cdot(f'\circ\phi)$ extends holomorphically also.  Furthermore,
the derivative of the extension does not vanish on the boundary.
Since $f'\circ\phi$ extends to the double of $\Ot$, and since $\Ot$
has real analytic boundary, it follows that
$f'\circ\phi$ extends holomorphically past the boundary of $\Ot$.
We conclude that $\phi'$ extends past the boundary of $\Ot$ at
all but the finitely many boundary points where $f`\circ\phi$ might
vanish, and at the vanishing points, $\phi$ maps the boundary of $\Ot$
to a cusp-like boundary point of $\O$.  Thus we conclude that $\O$
must have piecewise real analytic boundary.

Next, to see that $\O$ is a generalized quadrature domain,
let $z(t)$ parameterize one of the boundary curve segments of $\O$.
Since $\ln |f(z(t))|\equiv 1$, it follows by differentiating with
respect to $t$ that
$$\frac{f'(z(t))}{f(z(t))} z'(t)= -
\left(\,\overline{f'(z(t))}/\,\overline{f(z(t))}\,\right) \overline{z'(t)}.
\tag5.2$$
Let $f_b$ be an Ahlfors map such that $f$ and $f_b$ generate the
meromorphic functions on the double of $\O$.  Since $f'$ extends
to the double as a meromorphic functions, we know that
$f'=R(f,f_b)$ for some rational function $R$.  Now $f=1/\overline{f}$
on the boundary of $\O$ since $f$ maps the boundary into the unit
circle.  The same is true for $f_b$.  Hence, (5.2) yields that
$$z'(t)= -1/R\left(1/\,\overline{f(z(t))}\,,\,1/\overline{f_b(z(t))}\,\right)
\left(\,\overline{f'(z(t))}/\,\overline{f(z(t))^2}\right) \overline{z'(t)},
$$
i.e., that $dz=\overline{H(z)}d\bar z$ on $b\O$ where $H$ extends
meromorphically to $\O$.  Following Aharonov and Shapiro \cite{1},
Gustafsson shows in \cite{12, p.~223} that this condition is equivalent
to being a generalized quadrature domain.  This completes the proof.

\subhead 6. Quadrature domains with respect to arc length measure\endsubhead
An analogous theorem to Theorem~1.5 can be proved for smooth quadrature
domains with respect to boundary arc length measure.  Suppose
$\O$ is a bounded $n$-connected domain in the plane bounded by $n$
non-intersecting $C^\infty$ simple closed curves.  We say that $\O$
is a {\it quadrature domain with respect to arc length measure\/} if
there exist finitely many points
$\{w_j\}_{j=1}^N$ in the domain and non-negative integers $n_j$ such
that complex numbers $c_{jk}$ exist satisfying
$$\int_{b\O} f\ ds = \sum_{j=1}^N\sum_{k=0}^{n_j} c_{jk}
f^{(k)}(w_j)\tag6.1$$
for every function $f$ in the Hardy space $H^2(b\O)$ of holomorphic
functions on $\O$ with square integrable boundary values on $b\O$ with
respect to arc length measure $ds$ (see \cite{2} for basic facts
about $H^2(b\O)$).  The techniques used in the
previous sections can be adapted to replace the Runge theorems used
by Gustafsson in the proofs of his more general results by density
theorems in $A^\infty$ for the Szeg\H o kernel.  Indeed, we can follow
Gustafsson's argument in \cite{13, p.~76} to the letter, noting that
Gustafsson's function $h$ can be taken to be a complex linear
combination of functions of the form $S(z,b)$ where $b$ ranges over
an open subset of $\O$.  This is because identity (2.2) shows that
$S(z,b)^2dz=-\overline{L(z,b)^2} d\bar z$,
and hence, $h\sqrt{dz}$ is a ``half-order differential'' where
$h(z)=S(z,b)$.  Similar reasoning reveals the same thing about
complex linear combinations of such functions.  We may now
follow Gustafsson's argument, using the fact that the complex linear
span of $\{S(z,b)\,:\, b\in\O\}$ is dense in $A^\infty(\O)$ in
place of the Runge-type approximation theorem he uses.
Gustafsson's functions $f_j$ on page~77 can also be approximated in
$A^\infty(\O)$ by functions in this linear span.  In this way, we
may construct a function $h$ in $A^\infty(\O)$ such that $h^2$ has
a single valued antiderivative $g$ which is as close
to the identity map in $C^\infty(\Obar)$ as we desire.  This yields
a quadrature domain with respect to arc length measure $g(\O)$ which
is conformally equivalent to $\O$ and as $C^\infty$ close to $\O$ 
as we desire.

\enddocument